
\input amstex
\documentstyle{amsppt}
\def\L{{\Cal L}}
\def\M{{\Cal M}}
\def\O{{\Cal O}}
\def\p{{\Bbb P}}
\def\g{{\gamma}}
\def\s{{\sigma}}

\def\Q{{\Cal Q}}
\def\P{{\Cal P}}
\def\X{{\Cal X}}
\def\H{{\Cal H}}
\def\N{{\Cal N}}
\def\I{{\Cal I}}
\def\E{{\Cal E}}
\def\F{{\Cal F}}
\def\G{{\Cal G}}
\def\bN{{\Bbb N}}
\def\Z{{\Bbb Z}}
\def\Pn{{\p}^n}
\def\Pthree{{\p}^3}

\def\ext{\mathop{\Cal Ext}}

\def\ra{\rightarrow}

\topmatter
\title Integral Subschemes of Codimension Two \endtitle
\author Scott Nollet \endauthor
\address University of California at Riverside \endaddress
\email nollet\@math.ucr.edu \endemail
\abstract In this paper we study the problem of finding all integral
subschemes in a fixed even linkage class $\L$ of subschemes in $\Pn$ of pure
codimension two. To a subscheme $X \in \L$ we associate two invariants
$\theta_X, \eta_X$. When taken with the height $h_X$, each of these
invariants determines the location of $X$ in $\L$, thought of as a poset
under domination. In terms of these invariants, we find necessary conditions
for $X$ to be integral. The necessary conditions are almost sufficient in
the sense that if a subscheme dominates an integral subscheme and satisfies
the necessary conditions, then it can be deformed with constant cohomology
to an integral subscheme.
\endabstract
\date April 10, 1995 \enddate
\endtopmatter
\document
\heading Introduction \endheading
In \cite{4}, Gruson and Peskine show that an integral arithmetically
Cohen-Macaulay curve in $\Pthree$ has numerical character without
gaps. Conversely they show that any admissible sequence without gaps
arises as the numerical character of a smooth connected arithmetically
Cohen-Macaulay curve. Arithmetically Cohen-Macaulay curves form just one
of many even linkage classes of curves in $\Pthree$. In the present paper,
we obtain a similar result for even linkage classes of pure codimension two
subschemes in $\Pn$ which are not arithmetically Cohen-Macaulay. \par
In the first section we briefly review the linkage theory of codimension
two subschemes in $\Pn$. We recall the notions of $\E$ and
$\N$-type resolutions, the cone construction,
Rao's correspondence, and criteria for domination in an even linkage class.
We close the section with the theorem that the Lazarsfeld-Rao property
holds for any such even linkage class which does not consist of arithmetically
Cohen-Macaulay subschemes. \par
In the second section, we recall Martin-Deschamps and Perrin's notion
of admissible character and the partial ordering of domination on them.
We show that there is a bijection between dominations $\g \leq_h \s$ and
functions $\eta, \theta: \Z @>>> \bN$ with certain properties (proposition 2.6
and proposition 2.9). \par
Section three is devoted to adapting the invariants of section two to
the geometric situation. Specifically, it is shown that for a subscheme
$X$ in an even linkage class $\L$ of pure codimension two subschemes in
$\Pn$ which are not arithmetically Cohen-Macaulay, either of the invariants
$\eta_X$ or $\theta_X$ determine the class of all subschemes that can
be deformed with constant cohomology to $X$ through subschemes in $\L$.
\par
In section four, we give a preliminary result towards studying the integral
subschemes in an even linkage class. If $X$ is a subscheme, we
may define $s_0(X)$ to be the least degree of a hypersurface $S$ on which
$X$ lies, $s_1(X)$ to be the least degree hypersurface on which $X$ lies
which does not contain $S$, and $t_1(X)$ to be the least degree of a
hypersurface which meets $S$ properly. If $X$ is an integral subscheme,
then $X$ lies on an integral hypersurface of minimal degree, hence
$s_1(X)=t_1(X)$. The main result of this section (theorem 4.7) is a
characterization in terms of $\theta$ of which subschemes satisfy this
latter property. \par
In section five, we prove our main results about integral subschemes
in a fixed even linkage class. We give necessary conditions for a
subscheme $X$ to be integral in terms of the invariants $\theta_X$ and
$s_0(X)$ (theorem 5.8). Further, we show that if $X$ is integral,
$X \leq Y$, and $Y$ satisfies the necessary conditions, then $Y$
deforms with constant cohomology through subschemes in $\L$ to an
integral subscheme ${\overline Y}$ (theorem 5.11). To close out the
section, we give several examples of different kinds of behavior of
integral subschemes in an even linkage class. \par
While the results above are good for integral subschemes, the situation
is not so simple for smooth connected subschemes (example 5.15).
To obtain a result like theorem 5.11 for smooth connected subschemes
of codimension two in $\Pn$, more conditions are needed. Even with
extra conditions, typical constructions of smooth subschemes won't
work when $n \geq 6$. \par
The ideas of this paper originated in my PhD Thesis, where the case of
curves in $\Pthree$ was studied. I would like to thank Robin Hartshorne for a
careful reading of that thesis, as his questions and suggestions led to many
improvements.
\par
%
%
\heading $1.$ Linkage Theory in Codimension Two \endheading
   In this section we review the main results of linkage theory for
subschemes in $\Pn$ of codimension two. The main reference for this
section is \cite{19}, where various results of linkage theory for
Cohen-Macaulay subschemes are generalized to subschemes of pure
codimension. The main result of
importance here is the fact that the Lazarsfeld-Rao property holds
for even linkage classes of subschemes of pure codimension two in
$\Pn$. We assume that the reader is familiar with the notion of
simple linkage and the equivalence relation that it generates (see
\cite{8,$\S 4$}). \par
\definition{Definition 1.1}
A sheaf $\F$ on $\Pn$ is called {\it dissoci\'e} if
it is a direct sum of line bundles.
\enddefinition
\definition{Definition 1.2}
Let $V \subset \Pn$ be a subscheme of codimension $2$. An
{\it $\E$-type resolution for} $V$ is an exact sequence
$$0 \ra \E \ra \Q \ra \O (\ra \O_V \ra 0)$$
such that $\Q$ is dissoci\'e and $H^1_*(\E)=0$.
\enddefinition
\definition{Remark 1.3}
The condition that $H^1_*(\E)=0$ is equivalent to the the condition that
$H^0_*(\Q) \ra I_V$ is surjective. In particular, any subscheme $V$ of
codimension two has an $\E$-type resolution obtained by sheafifying a free
graded surjection onto the ideal $I_V$. The theorem of Auslander and Buchsbaum
shows that $V$ is locally Cohen-Macaulay $\iff \E$ is locally free. Further,
it can be shown that
$V$ is of pure codimension two $\iff \ext^1(\E^\vee, \O)=0$. This second
equivalence generalizes to subvarieties of higher codimension
\cite{19, corollary 1.18}.
\enddefinition
\definition{Definition 1.4}
If $V \subset \Pn$ is a subscheme of codimension $2$, then an
{\it $\N$-type resolution for $V$} is an exact sequence
$$0 \ra \P \ra \N \ra \O (\ra \O_V \ra 0)$$
where $\P$ is dissoci\'e, $\N$ is reflexive, $\ext^1(\N, \O)=0$ and
$H^1_*(\N^\vee)=0$.
\enddefinition
\definition{Remark 1.5}
Unlike the situation for $\E$-type resolutions, not every subscheme of
codimension two has an $\N$-type resolution. In fact, a subscheme
$V \subset \Pn$ of codimension $\geq 2$ has an $\N$-type resolution if and only
if $V$ is of pure codimension $2$. A similar statement holds in higher
codimension \cite{19, corollary 1.20}.
\enddefinition
\proclaim{Proposition 1.6}
Let $V \subset \Pn$ be a subscheme of pure codimension $2$ which is contained
in a complete intersection $X$ of two hypersurfaces with degrees summing to
$d$.
Let
$$0 \ra \G_2 \ra \G_1 \ra \G_0 (\ra \O_V \ra 0)$$
be an $\E$-type (respectively $\N$-type) resolution for $V$ and let $\F_.$ be a
Koszul resolution for $X$. Then there is a morphism of complexes
$\alpha:\F_. @>>> \G_.$ induced by the inclusion $\I_X @>>> \I_V$. The
mapping cone of the morphism $\alpha^\vee(-d)$ gives a resolution for the
subscheme $W \subset \Pn$ which is linked to $V$ by $X$. If the induced
isomorphism $\G_0^\vee(-d) @>>> \F_0^\vee(-d)$ is split off, the resulting
resolution is an $\N$-type (respectively $\E$-type) resolution for $W$.
\endproclaim
\demo{Proof}
This follows from \cite{19, $\S 1$}.
\enddemo
\definition{Remark 1.7} Let $0 \ra \E_V \ra \Q \ra \O$ be an $\E$-type
resolution
for a subscheme $V$ as in the proposition. If $V$ is linked to $W$ and $W$ is
linked to $Z$, then applying proposition 1.6 twice produces an
$\E$-type resolution for $Z$ such that $\E_Z = \E_V(h) \oplus \F$ where
$h \in \Z$ and $\F$ is dissoci\'e. This motivates the definition of stable
equivalence and the theorems which follow.
\enddefinition
\definition{Definition 1.8} Two reflexive sheaves $\E_1$ and $\E_2$ on $\Pn$
are {\it stably equivalent} if there exist dissoci\'e sheaves $\Q_1,\Q_2$
and $h \in \Z$ such that $\E_1 \oplus \Q_1 \cong \E_2(h) \oplus \Q_2$.
This is an equivalence relation among reflexive sheaves on $\Pn$.
\enddefinition
\proclaim{Theorem 1.9}
There is a bijective correspondence between even linkage classes of
purely two-codimensional subschemes of $\Pn$ and stable equivalence classes
of reflexive sheaves $\E$ on $\Pn$ such that $H^1_*(\E)=0$ and
$\ext^1(\E^\vee, \O)=0$ via $\E$-type resolutions.
\endproclaim
\demo{Proof} This is \cite{19, theorem 2.11}, where the ideas of Rao's
original proof \cite{22} are extended to the case of subschemes which
are not locally Cohen-Macaulay.
\enddemo
\proclaim{Theorem 1.10}
There is a bijective correspondence between even linkage classes of purely
two-codimensional subschemes of $\Pn$ and stable equivalence classes of
reflexive sheaves $\N$ on $\Pn$ such that $H^1_*(\N^\vee)=0$ and
$\ext^1(\N, \O)=0$ via $\N$-type resolutions.
\endproclaim
\demo{Proof} This follows from theorem 1.9 and proposition 1.6 once we note
that
two reflexive sheaves are stably equivalent if and only if their duals are
stably equivalent. \enddemo
\proclaim{Proposition 1.11}
Let $\Omega$ be a stable equivalence class of reflexive sheaves on $\Pn$.
Let $\F_0$ be a sheaf in $\Omega$ of minimal rank. Then for each
$\G \in \Omega$, there exists $h \in \Z$ and a dissoci\'e sheaf $\Q$ such
that $\G \cong \F_0(h) \oplus \Q$. In particular, the isomorphism class of
$\F_0$ is uniquely determined up to twist.
\endproclaim
\demo{Proof} This is \cite{19, proposition 2.3}. \enddemo
\definition{Notation 1.12}
Let $\L$ be an even linkage class of codimension two subschemes in
$\Pn$. If $\L$ is associated with the stable equivalence class $\Omega$
via $\E$-type resolution (as in theorem 1.9) and $\E_0$ is a sheaf
of minimal rank for $\Omega$ (as in proposition 1.11), we will say that
$\L$ {\it corresponds to $[\E_0]$ via $\E$-type resolution}.
Similarly, if $\L$ is associated to $\Omega$ via $\N$-type resolution
and $\N_0$ is a sheaf of minimal rank for $\Omega$, we will say that
$\L$ {\it corresponds to $[\N_0]$ via $\N$-type resolution}.
\enddefinition
The above theorems give a classification of even linkage classes. In the
remainder of this section, we describe the structure of an even linkage
class.
\definition{Definition 1.13}
Let $V \subset \Pn$ be a subscheme of pure codimension two.
Let $S$ be a hypersurface of degree $s$ which contains $V$ and let $h$
be an integer. We say that $W$ is obtained from $V$ by an {\it elementary
double
link of height} $h$ {\it on} $S$ if there are hypersurfaces $T_1, T_2$ of
respective degrees $t_1, t_1+h$ which meet $S$ properly in such a way that
$S \cap T_1$ links $V$ to a subscheme $Z$, which in turn is linked to
$W$ by $S \cap T_2$. In the special case when $h \geq 0$ and $T_2=T_1 \cup H$
where $H$ is a hypersurface of degree $h$, we say that $Y$ is obtained from
$X$ by a {\it basic double link of height} $h$ {\it on} $S$. In either case,
we say that the double link has type $(s,h)$.
\enddefinition
\proclaim{Proposition 1.14}
Suppose that $V \subset \Pn$ is a subscheme of pure codimension two and
that $W$ is obtained from $X$ by an elementary double link of type $(s,h)$.
If $V$ has an $\E$-type (respectively $\N$-type) resolution of the form
$$0 \ra \G_2 \ra \G_1 \ra \O (\ra \O_V \ra 0),$$
then $W$ has an $\E$-type (respectively $\N$-type) resolution of the
form
$$0 \ra \G_2(-h) \oplus \O(-s-h) \ra \G_1(-h) \oplus \O(-s) \ra \O (\ra \O_W
\ra 0). $$
\endproclaim
\demo{Proof} Applying proposition 1.6 twice, we obtain such a resolution
with the extra summand $\O(-t_1-h)$ appearing in the first two terms.
Because the same hypersurface $S$ was used for both links, we can take
the induced map between these summands to be the identity. Splitting off
this summand, we are left with the resolution stated.
\enddemo
\definition{Definition 1.15}
Let $X,X^\prime$ be subschemes of pure codimension two
in $\Pn$. We say that $X^\prime$ {\it dominates} $X$ {\it at height}
$h \geq 0$ if $X^\prime$ can be obtained from $X$ by a sequence of basic
double links with heights summing to $h$, followed by a deformation which
preserves cohomology and even linkage class. In this case we write
$X \leq_h X^\prime$, or simply $X \leq X^\prime$ if $h$ is not specified.
\enddefinition
\definition{Definition 1.16}
If $f:\Z \ra {\bN}$ is a function such that $f(n)=0$ for
$n << 0$, then we define the function $f^\#$ by $f^\#(a)=\sum_{n \leq a}f(n)$.
\enddefinition
\proclaim{Proposition 1.17}
Let $X,Y \subset \Pn$ be subschemes of pure codimension two and let $h \geq 0$
be an integer. Then $X \leq_h Y$
if and only if there exist $\N$-type resolutions
$$0 \ra \oplus\O(-m)^{r(m)} \ra \N \ra \I_X(h_1) \ra 0$$
$$0 \ra \oplus\O(-m)^{s(m)} \ra \N \ra \I_Y(h_2) \ra 0$$
such that $h_2-h_1=h$ and $r^\# \geq s^\#$.
\endproclaim
\demo{Proof} This is \cite{19, proposition 3.9}.
\enddemo
\proclaim{Corollary 1.18}
if $X \leq_h Y$ and $Y \leq_k Z$, then $X \leq_{h+k} Z$. In particular,
domination is transitive.
\endproclaim
\demo{Proof} This is \cite{19, corollary 3.10}
\enddemo
\definition{Definition 1.19}
Let $\L$ be an even linkage class of subschemes of pure
codimension two in $\Pn$. We say that $\L$ has the
{\it Lazarsfeld-Rao property} if $\L$ has
a minimal element $X_0$ such that $X_0 \leq X$ for each $X \in \L$.
\enddefinition
\proclaim{Theorem 1.20}
Let $\L$ be an even linkage class of codimension two subschemes of $\Pn$
corresponding to the stable equivalence class $[\N_0]$ via $\N$-type
resolutions, with $\N_0 \neq 0$. Then there exists a function $q:\Z \ra \bN$
and a subscheme $X_0 \in \L$ such that $X_0 \leq Y$ for all $Y \in \L$ and
$X_0$ has an $\N$-type resolution of the form
$$0 \ra \oplus\O(-n)^{q(n)} \ra \N_0 \ra \O_{\Pn}(h) (\ra \O_{X_0}(h) \ra 0).$$
In particular, $\L$ has the Lazarsfeld-Rao property. $X_0^\prime$ is
another such minimal element of $\L \iff X_0^\prime$ has an $\N$-type
resolution of the same form as that of $X_0 \iff X_0$ deforms to $X_0^\prime$
with constant cohomology through subschemes in $\L$.
\endproclaim
\demo{Proof} This is \cite{19, theorem 3.26}. \enddemo
\proclaim{Theorem 1.21}
Let $\L$ be an even linkage class corresponding to the stable
equivalence class $[\N_0]$ via $\N$-type resolution, with
$\N_0 \neq 0$. Let $X_0$ be a minimal subscheme for $\L$. If $S$ is a
hypersurface of minimal degree which contains $X_0$, then there exists
a hypersurface $T$ containing $X_0$ such that $S \cap T$ links $X_0$ to
a subscheme $Y_0$ which is minimal for its even linkage class.
\endproclaim
\demo{Proof} This is a simplified version of \cite{19, theorem 3.31}.
\enddemo
%
%
\heading $2.$ Domination of Admissible Characters \endheading
In this section we recall Martin-Deschamps and Perrin's notion of
admissible characters, and the partial ordering of domination on
them \cite{13,V}. In this section we give properties of the partial
ordering along with alternative ways of describing it. The relationship
between domination of admissible characters and domination of subschemes
in $\Pn$ will be explained in the next section. \par
\definition{Definition 2.1} A {\it character} is a function $\s:\Z @>>> \Z$
which has finite support and satisfies $\sum_{l \in \Z}{\s(l)}=0$. A
character $\s$ is said to be {\it admissible} if it satisfies
\roster
\item $\s(l)=0$ for $l < 0$
\item $\s(0)=-1$
\item Setting $s_0(\s)=inf\{l \geq 0:\s(l) \neq -1\}$, we have $\s(s_0(\s))
\geq 0$
\item Setting $s_1(\s)=inf\{l \geq s_0(\s):\s(l) \neq 0\}$, we have
$\s(s_1(\s)) > 0$
\endroster
\enddefinition
\definition{Definition 2.2} Let $\g,\g^\prime$ be two admissible
characters. If $h \geq 0$ is an integer, we say that {\it $\g^\prime$
dominates $\g$ at height $h$} if
\roster
\item $s_0(\g) \leq s_0(\g^\prime) \leq s_0(\g)+h$
\item $\g^\prime(l) \geq 0$ for $s_0(\g^\prime) \leq l < s_0(\g)+h$
\item $\g^\prime(l) \geq \g(l-h)$ for $l \geq s_0(\g)+h$
\endroster
In this case we write $\g \leq_h \g^\prime$. If $h$ isn't specified,
we simply write $\g \leq \g^\prime$.
\enddefinition
\proclaim{Proposition 2.3} If $\g,\s,\tau$ are three admissible
characters and $\g \leq_h \s$ and $\s \leq_k \tau$, then
$\g \leq_{h+k} \tau$. In particular, domination gives a partial
order relation on the set of admissible characters.
\endproclaim
\demo{Proof} This is \cite{13, V, proposition 2.2}. \enddemo
\definition{Definition 2.4} Let $f:\Z @>>> \bN$ be a function. then
$f_o=max\{l:f(l)>0\}$ if such an $l$ exists, otherwise $f_o=-\infty$. Similarly
we define $f_a=min\{l:f(l)>0\}$ if such an $l$ exists, otherwise $f_a=\infty$.
\enddefinition
\definition{Definition 2.5} Let $f:\Z @>>> \bN$ be a function.
We say that $f$ is {\it connected in degrees $\geq a$ (respectively $>a$)}
if $f(b)>0$ for some $b > a$ implies that $f(l)>0$ for all $a \leq l \leq b$
(respectively $a < l \leq b$). Similarly we say that $f$ is
{\it connected in degrees $\leq b$ (respectively $<b$)} if $f(a)>0$
for some $a < b$ implies that $f(l)>0$ for $a \leq l \leq b$ (respectively
$a \leq l < b$). $f$ is said to be {\it connected about an interval
$[a,b]$} if $f$ is connected in degrees $\leq b$, $f$ is connected in
degrees $\geq a$ and $f$ is nonzero on $[a,b]$. $f$ is {\it connected} if
it is connected about the interval $[f_a,f_o]$.
\enddefinition
\proclaim{Proposition 2.6} Let $\g$ be an admissible character,
$h \geq 0$ be an integer and $\s:\Z @>>> \Z$ be a function.
Then $\s$ is an admissible character such that $\s \geq_h \g$ if and only
if the function $\eta$ given by
$$\eta(l)=\s(l)-\g(l-h)+{l \choose 0}-{{l-h} \choose 0}$$
is nonnegative, is connected in degrees $< s_0(\g)+h$ and satisfies
$\sum{\eta(l)}=h$. This function is denoted $\eta_{\g, \s, h}$ or simply
$\eta$.
\endproclaim
\demo{Proof}
Assume that $\s \geq_h \g$ is an admissible character and define $\eta$
as above. A simple calculation shows that $-\g(l-h)+{l \choose 0}-
{{l-h} \choose 0}=1$ for $0 \leq l < s_0(\g)+h$. It follows that
$\eta(l)=0$ for $0 \leq l < s_0(\s)$ and that $\eta(l)=\s(l)+1$ for
$s_0(\s) \leq l < s_0(\g)+h$. Since $\s(l) \geq 0$ in this last range
by definition of domination, we see that $\eta$ is connected
in degrees $< s_0(\g)+h$. For $l \geq s_0(\g)+h$, we have that
$\eta(l)=\s(l)-\g(l-h)$, which is $\geq 0$ for such $l$ by definition
of domination. It follows that $\eta(l) \geq 0$ for all $l \in \Z$.
Finally, $\sum{{l \choose 0}-{{l-h} \choose 0}}=h$, so the
fact that $\g$ and $\s$ are characters shows that $\sum{\eta(l)}=h$.
\par
Conversely, suppose that $\eta:\Z \ra \bN$ satisfies the conditions of the
proposition. If we define $\s$ by the formula for $\eta$, the fact that
$\sum{\eta(l)}=h$ shows that $\s$ is a character. Since $\eta$ sums to $h$
and is connected in degrees $< s_0(\g)+h$, we see that $\eta(l)=0$ for
$l < s_0(\g)$. It follows that $\s(l)=0$ for $l<0$ and $\s(l)=-1$ for
$0 \leq l < s_0(\g)$. In particular, $\s$ satisfies the first two conditions
for admissibility and $s_0(\s) \geq s_0(\g)$. Since $\s(s_0(\g)+h) \geq 0$,
we also see that $s_0(\s) \leq s_0(\g)+h$, hence the first condition for
domination holds. The third condition of domination holds because
$\s(l)-\g(l-h)=\eta(l) \geq 0$ for $l \geq s_0(\g)+h$. \par
To check the second domination condition, we consider two cases.
If $\eta(l)=0$ for $0 \leq l < s_0(\g)+h$, then
$\s(l)=-1$ for $l$ in this range, so $s_0(\s)=s_0(\g)+h$ and the second
domination condition holds vacuously. If this is not the case, let
$N = \text{min}\{l:\eta(l)>0\}$. The connectedness condition on $\eta$
shows that $\eta(l)>0$ for $N \leq l < s_0(\g)+h$, which implies that
$s_0(\s)=N$ and $\s(l) \geq 0$ for $s_0(\s) \leq l < s_0(\g)+h$. This checks
the second condition, so all three conditions for domination hold. \par
Finally, we check the last two conditions for admissibility. From the
domination conditions, we see that $\s(l) \geq 0$ for
$s_0(\s) \leq l < s_0(\g)+h$, which immediately gives the third
admissibility condition. Further, we have that $\s(l) \geq 0$ for
$s_0(\g)+h \leq l \leq s_1(\g)+h$ and that $\s(s_1(\g)+h)>0$. It follows
that the last admissibility condition holds, finishing the proof.
\enddemo
\proclaim{Proposition 2.7}
Let $\g \leq_h \tau$ and $\g \leq_k \s$ be two dominations of admissible
characters. Then $\tau \leq_{k-h} \s$ if and only if
$\eta_{\g, \s}(l)-\eta_{\g, \tau}(l-k+h) \geq 0$ for all $l \in \Z$,
in which case we have the formula
$$\eta_{\tau, \s, k-h}(l)=\eta_{\g, \s, k}(l)-\eta_{\g, \tau, h}(l-k+h).$$
\endproclaim
\demo{Proof}
If $k \geq h$, then a simple calculation shows gives the formula above.
In the case that $\tau \leq_{k-h} \s$, we have that $\eta_{\tau, \s, k-h}$ is
nonnegative by proposition 2.6, so the forward implication is clear. \par
Conversely, suppose that the function $\eta=\eta_{\tau, \s, k-h}$ as defined by
the formula is nonnegative. It is clear that $\sum{\eta(l)}=k-h$, so we have
that $k-h \geq 0$. For $l < s_0(\tau)+k-h$ we have that $\eta_{\g,\tau}(l)=0$,
so $\eta(l)=\eta_{\tau,\s}(l)$ in this range. Since $\eta_{\tau,\s}$ is
connected in degrees $< s_0(\tau)+k-h$ (in fact, in degrees $< s_0(\g)+k$),
this holds for $\eta$ as well. We deduce the domination $\tau \leq_{k-h} \s$
from proposition 2.6.
\enddemo
\definition{Remark 2.8} In the the proof of proposition 2.6, we showed that
$\eta(l)=\s(l)+1$ for $s_0(\s) \leq l < s_0(\g)+h$. If we had used
simply $\s(l)$ instead, we would get a second function $\theta$. It turns
out that while $\eta$ is
easy to deal with algebraically, the function $\theta$ has some nice
geometric properties, as we will see in the next section. For $\theta$
we have the following analogous proposition.
\enddefinition
\proclaim{Proposition 2.9} Let $\g$ be an admissible character and
$h \geq 0$ an integer. Then there is a bijection between admissible
characters $\s \geq \g$ and functions $\theta:\Z @>>> \bN$ such that
$\sum{\theta(l)}=m \leq h$ and $\theta(l)=0$ for $l < s_0(\g)+m$.
If $\s \geq_h \g$ is given and $\eta$ is the corresponding function from
proposition 2.6, then $\theta$ is given by
$$\theta(l)=\eta(l)-{{l-s_0(\s)} \choose 0}+{{l-s_0(\g)-h} \choose 0}.$$
This function is denoted $\theta_{\g, \s, h}$, or just $\theta$.
\endproclaim
\demo{Proof}
Let $\s \geq_h \g$ be the given domination of admissible characters and
let $\eta,\theta$ be the two functions defined via propositions 2.6 and
proposition 2.9. As noted in remark 2.8, we have that $\theta(l)=\s(l)$ for
$s_0(\s) \leq l < s_0(\g)+h$.
Since $\s(l) \geq 0$ for $l$ in this range, we conclude that
$\theta(l) \geq 0$ for all $l \in \Z$. Clearly $m=\sum{\theta(l)} \leq
\sum{\eta(l)}=h$. Since $\eta(l)=0$ for $l < s_0(\s)=s_0(\g)+m$, this
is also true for $\theta$. \par
   Conversely, suppose that we are given $\theta$ as in the proposition, with
$m = \sum{\theta(l)}$. To give the corresponding $\s \geq_h \g$, it
suffices to produce $\eta$ with the conditions of proposition 2.6.
For $l \in \Z$, we define
$$\eta(l)=\theta(l)+{{l-s_0(\g)-m} \choose 0}-{{l-s_0(\g)-h} \choose 0}.$$
It is easily seen that $\eta(l) \geq 0$ for each $l \in \Z$ and that
$\sum{\eta(l)}=h$. Because $\theta(l)=0$ for $l < s_0(\g)+m$, the formula
shows that $\eta$ also vanishes in this range. The formula also
shows that $\eta(l) \geq 1$ for $s_0(\g)+m \leq l < s_0(\g)+h$, hence
$\eta$ is connected in degrees $< s_0(\g)+h$.
\enddemo
\proclaim{Proposition 2.10}
Let $\g \leq_h \tau$ and $\g \leq_k \s$ be two dominations of admissible
characters. Then $\tau \leq_{k-h} \s$ if and only if
$\theta(l)=\theta_{\g, \s}(l)-\theta_{\g, \tau}(l-k+h) \geq 0$ for all
$l \in \Z$ and $\sum{\theta(l)} \leq k-h$, in which case we have
$\theta_{\tau, \s, k-h}=\theta.$
\endproclaim
\demo{Proof}
First assume that $\tau \leq_{k-h} \s$. In this case the definition of
domination gives that $h \leq k$ and
$s_0(\tau) \leq s_0(\s) \leq s_0(\tau)+k-h$.
Since $s_0(\s)=s_0(\g)+\sum{\theta_{\g, \s}(l)}$ and
$s_0(\tau)=s_0(\g)+\sum{\theta_{\g, \tau}}$, the second inequality is
equivalent to $0 \leq \sum{\theta(l)} \leq k-h$. It is an easy calculation
that $\theta_{\tau, \s}=\theta$. Since $\theta_{\tau, \s}$ is
nonnegative, we deduce the forward direction. \par
Conversely, assume the two conditions. As noted in the preceding paragraph,
$\sum{\theta(l)} \leq k-h$ is equivalent to $s_0(\s) \leq s_0(\tau)+k-h$.
A tedious calculation gives the formula
$$\eta_{\tau, \s}(l)=\theta(l)+{{l-s_0(\s)} \choose 0}-{{l-s_0(\tau)-k+h}
\choose 0},$$
which comes as no surprise given the statement of proposition 2.8. Now the
fact that $\theta$ is nonnegative and the inequality
$s_0(\s) \leq s_0(\tau)+k-h$ show that $\eta_{\tau, \s}$ is nonnegative.
Applying proposition 2.10 gives that $\tau \leq_{k-h} \s$.
\enddemo
%
%
\heading $3.$ Geometric Invariants \endheading
Propositions 2.6 and 2.9 give two descriptions of the difference
between a larger and smaller admissible character. Now we apply
this to the geometric situation. In what follows, we will use
$n$th difference functions quite a bit, so we make the following
abbreviation to simplify matters. \par
\definition{Notation 3.1}
Let $\F$ be a coherent sheaf on $\Pn$. We will use the
abbreviation $\Delta^m\F$ to denote the $m$th difference function of
the cohomology function $h^0(\F(l))$. In other words,
$\Delta^m\F(l)=\Delta^m(h^0(\F(l)))$.
\enddefinition
\proclaim{Proposition 3.2} Let $X \subset \Pn$ be a subscheme of
codimension $\geq 2$. Then the function $\g:\Z @>>> \Z$ defined by
$$\g(l)=\Delta^n\I_X(l) - {l \choose 0}$$
is an admissible character.
\endproclaim
\demo{Proof}
Consider the function $\varphi(l)=h^0(\I_X(l))-h^0(\O_{\Pn}(l))$.
Clearly $\varphi(l)=0$ for $l<0$. By Serre's vanishing theorem,
we see that $\varphi(l)=-P_X(l)$ for $l >> 0$, where $P_X$ is the
Hilbert polynomial for $X$. Since $X$ is of codimension $\geq 2$,
the Hilbert polynomial has degree $\leq n-2$, hence we have that
$\Delta^{(n-1)}(P_X(l))=0$ for all $l$. It follows that the function
$\Delta^{(n-1)}\varphi(l)$ has finite support. To see that $\g$ is
a character, we note more generally that the first difference of a
function of finite support is a character. \par
To check admissibility, we first note that $\g(l)=0$ for $l<0$ because
this holds for $\varphi$. If we set
$s_0=min\{l:h^0(\I_X(l)) \neq 0\}$, we see that $\g(l)=-1$ for
$0 \leq l < s_0$ because $\Delta^n\O_{\Pn}(l)= {l \choose 0}$. Choosing
$0 \neq v \in H^0(\I_X(s_0))$ gives rise to an injection
$\tau:\O(-s_0) \hookrightarrow \I_X$, whose image consists of all
multiples of an equation for a hypersurface $S$ of minimal degree which
contains $X$. Letting $s_1$ be the least degree of a hypersurface
$T$ which contains $X$ but is not a multiple of $S$ (such hypersurfaces
exist because $X$ has codimension $\geq 2$), we see that $s_1$ is the
least twist where $H^0(\tau)$ is not surjective. Since
$\Delta^n\O(-s)(l)={{l-s} \choose 0}$, we have that $\g(l)=0$ for
$s_0 \leq l < s_1$ and that $\g(s_1)>0$. This shows that $\g$ is
an admissible character.
\enddemo
\definition{Definition 3.3} Let $X \subset \Pn$ be a subscheme of
pure codimension two. The admissible character
$\g(l)=\Delta^n\I_X(l)-{l \choose 0}$ of
proposition 3.2 is called the {\it $\g$-character} of $X$, and is
denoted $\g_X$. We define $s_0(X)=s_0(\g_X)$, $s_1(X)=s_1(\g_X)$ and
$e(X)=max\{l:H^{n-2}(\O_{X_0}(l)) \neq 0\}$. $t_1(X)$ is
the least degree of a hypersurface $T \supset X$ such that there is
a hypersurface $S \supset X$ of degree $s_0(X)$ which meets $T$ properly.
The {\it higher Rao modules of $X$} are the defined by
$M_X^i=H^i_*(\I_X)$ for $1 \leq i \leq n-2$. These are graded modules
over the homogeneous coordinate ring $S=H^0_*(\O_{\Pn})$ of $\Pn$.
\enddefinition
\definition{Remark 3.4}
The geometric meaning of the integers $s_0(X)$ and $s_1(X)$ is described
in the proof of proposition 3.2. We always have the inequalities
$s_0(X) \leq s_1(X) \leq t_1(X)$.
\enddefinition
\definition{Remark 3.5} Let $X$ and $Y$ be two subschemes of $\Pn$ having
codimension $\geq 2$. In this case, the condition that $\g_X \leq_h \g_Y$ has
a simple formulation in terms of the ideal sheaves. The function $\eta$ of
proposition 2.6 can be written
$$\eta(l)=\g_Y(l)-\g_X(l-h)+{l \choose 0}-{{l-h} \choose 0}$$
$$=\Delta^n\I_Y(l)-\Delta^n\I_X(l-h).$$
Applying proposition 2.6, it follows that $\g_X \leq_h \g_Y$ if and only
if the function $\Delta^n\I_Y(l)-\Delta^n\I_X(l-h)$ is nonnegative,
connected in degrees $< s_0(\g_X)+h$ and sums to $h$.
\enddefinition
\proclaim{Proposition 3.6} Let $X \subset \Pn$ be a subscheme and suppose
that $Y$ is obtained from $X$ by an elementary double link of type
$(s,h)$. Then $\g_X \leq_h \g_Y$, $M_Y^i \cong M_X^i(-h)$ for
each $1 \leq i \leq n-2$ and $\eta=\eta_{\g_X, \g_Y, h}$ is given by
$$\eta(l)={{l-s} \choose 0}-{{l-s-h} \choose 0}.$$
\endproclaim
\demo{Proof}
   Suppose that $Y$ is obtained from $X$ by an elementary double
link of type $(s,h)$. $X$ has an $\N$-type resolution that gives
rise to the exact sequence
$$0 \ra \P \ra \N \ra \I_X \ra 0.$$
By proposition 1.14, there is an $\N$-type resolution for $Y$ that
gives an exact sequence
$$0 \ra \P(-h) \oplus \O(-s-h) \ra \N(-h) \oplus \O(-s) \ra \I_Y \ra 0.$$
Twisting the second sequence by $h$ and adding the trivial summand
$\O(-s+h)$ to the first resolution gives a pair of exact sequences
$$0 \ra \P \oplus \O(-s+h) \ra \N \oplus \O(-s+h) \ra \I_X \ra 0$$
$$0 \ra \P \oplus \O(-s) \ra \N \oplus \O(-s+h) \ra \I_Y(h) \ra 0.$$
{}From these sequences we see that $M_Y^i \cong H^i_*(\N(-h)) \cong M_X^i(-h)$
for $1 \leq i \leq n-2$. Because the exact sequences above are also exact
on global sections, we have that
$\Delta^n\I_Y(l)=\Delta^n\I_X(l-h)+{{l-s} \choose 0}- {{l-s-h} \choose 0}$,
hence $\eta(l)={{l-s} \choose 0}-{{l-s-h} \choose 0}$ by remark 3.5.
This function is clearly nonnegative and sums to $h$. The fact that $\eta$ is
connected in degrees $\leq s_0(X)+h$ follows from the fact
that $s \geq s_0(X)$. It follows from proposition 2.6 that $\g_X \leq_h \g_Y$.
\enddemo
\proclaim{Corollary 3.7} With the hypotheses of Proposition 3.6, further
assume that $h=1$. Then $\theta=\theta_{\g_X, \g_Y, 1}$ is given by the
following:
\roster
\item If $s=s_0(X)$, then $\theta(l)=0$ for all $l \in \Z$
\item If $s > s_0(X)$, then $\theta(l)={{l-s} \choose 0}-{{l-s-1} \choose 0}$
for all $l \in \Z$
\endroster
\endproclaim
\demo{Proof}
This is a simple consequence of proposition 3.6 and the formula
relating $\theta$ and $\eta$, once we note that $s_0(Y) > s_0(X) \iff
s > s_0(X)$.
\enddemo
\proclaim{Theorem 3.8} Let $\L$ be an even linkage class of codimension
two subschemes of $\Pn$. Let $X,Y \in \L$ and let $l \geq 0$ be an
integer. Then $X \leq_h Y$ if and only if $\g_X \leq_h \g_Y$ and
$M_Y^i \cong M_X^i(-h)$ for each $1 \leq i \leq n-2$.
\endproclaim
\demo{Proof}
    If $Y$ is obtained from $X$ by a basic double link of type $(s,h)$,
then $\g_X \leq_h \g_Y$ and $M_Y^i \cong M_X^i(-h)$ by
proposition 3.6. Since both of these conditions are transitive, it
is clear that if $Y$ is obtained from $X$ by a sequence of basic double links
with heights summing to $h$, then $\g_X \leq_h \g_Y$ and $M_Y^i \cong
M_X^i(-h)$. These invariants are fixed under a deformation which
preserves cohomology and even linkage class, so we have proved
the forward direction. \par
   Conversely, suppose that $\g_X \leq_h \g_Y$ and $M_Y^i \cong M_X^i(-h)$ for
$1 \leq i \leq n-2$. Let $\N_0$ be a reflexive sheaf such that $\L$ corresponds
to $[\N_0]$ via $\N$-type resolutions. Since $X$ and $Y$ are in $\L$, we
can use Rao's correspondence (theorem 1.10) to find $\N$-type resolutions
that give exact sequences
$$0 \ra P_X \ra Q_X \oplus \N_0 \ra \I_X \ra 0$$
$$0 \ra P_Y \ra Q_Y \oplus \N_0 \ra \I_Y(k) \ra 0$$
for some $k \in \Z$. If $\N_0 \neq 0$, then some of the higher Rao
modules of $X$ are nonzero, and the condition $M_Y^i \cong M_X^i(-h)$
shows that $k=h$. If all the higher Rao modules are zero, then $\N_0=0$
and we can twist the dissoci\'e sheaves $P_Y$ and $Q_Y$ by $h-k$ to
assume $k=h$. In either case, we can find such resolutions with $h=k$.
   \par
Adding $Q_Y$ to the first sequence and $Q_X$ to the second,
we obtain exact sequences
$$0 \ra \oplus\O(-l)^{a(l)} \ra \N \ra \I_X \ra 0$$
$$0 \ra \oplus\O(-l)^{b(l)} \ra \N \ra \I_Y(h) \ra 0.$$
The condition $\g_X \leq_h \g_Y$ tells us that
$\Delta^n\I_Y(l) \geq \Delta^n\I_X(l-h)$ for all $l \in \Z$ by remark 3.5.
Noting that the two sequences above
are exact on global sections, it follows that
$\Delta^n(\oplus\O(-l)^{a(l)}) \geq \Delta^n(\oplus\O(-l)^{b(l)})$.
Since $\Delta^n\O(-l)={l \choose 0}$, this condition is equivalent
to the condition that $a^\#(l) \geq b^\#(l)$ for all $l \in \Z$.
By proposition 1.17, we conclude that $X \leq_h Y$, finishing the
proof.
\enddemo
\proclaim{Proposition 3.9} Let $X \subset \Pn$ be a subscheme of
pure codimension two having $\g$-character $\g_X$ and let $h \geq 0$ be
an integer. If $\s$ is an admissible character with $\g_X \leq_h \s$, then
there exists a subscheme $Y \subset \Pn$ which is obtained from $X$ by
a sequence of basic double links of height one such that
$X \leq_h Y$ and $\g_Y=\s$.
\endproclaim
\demo{Proof}
We induct on $h$. The base case is easy, since it is easy to check that
if $\g_X \leq_0 \s$ is a domination of admissible characters, then
in fact $\g_X = \s$. It follows that we may take $Y=X$ in this case.
\par
We now assume that $h>0$. Set $\eta=\eta_{\g_X,\s}$.
Since $\eta$ is connected in degrees $< s_0(X)+h$ by proposition 2.6, we
have that $\eta_o-h+1 \geq s_0(X)$, so we can obtain a subscheme $Z$ from
$X$ by a basic double link of type $(s,1)$, where $s=\eta_o-h+1$.
{}From propositions 2.7 and 3.6 we see that $\g_Z \leq_{h-1} \s$. By induction
hypothesis we can obtain a subscheme $Y \geq_{h-1} Z$ which is obtained
from $Z$ by a sequence of basic double links and satisfies $\g_Y = \s$.
Extending the chain of basic double links by the link from $X$ to $Z$
completes the proof.
\enddemo
\proclaim{Corollary 3.10} Let $\L$ be an even linkage class of subschemes
of pure codimension two in $\Pn$ and let $X \in \L$. Then there is a bijection
between cohomology preserving deformation classes of subschemes $Y \in \L$
such that $X \leq_h Y$ and admissible characters $\s \geq_h \g_X$.
\endproclaim
\demo{Proof}
If $\{Y_h\}_{h \in H}$ is such a deformation class, then $\g_{Y_h}$ is
a constant admissible character $\geq_h \g_X$. On the other hand, if
$\s \geq_h \g_X$, then we can apply proposition 3.9 to obtain
a subscheme $Y \geq_h X$ with $\g_Y=\s$. If $Y^\prime$ is another such
subscheme, then $\g_Y=\g_{Y^\prime}$ and
$M_{Y^\prime}^i \cong M_X^i(-h) \cong M_Y^i$ (by theorem 3.8). Applying
theorem 3.8 again, we have that $Y \leq_0 Y^\prime$, which is equivalent
to saying that there is a deformation from $Y$ to $Y^\prime$ through
subschemes in $\L$ with constant cohomology.
\enddemo
\definition{Remark 3.11} In \cite{2, theorem 3.5.7}, Bolondi and Migliore
prove essentially the same result in the case that $X$ is a minimal
element for a non-ACM even linkage class of locally Cohen-Macaulay
subschemes of codimension two in $\Pn$.
\enddefinition
Let $\L$ be an even linkage class of pure codimension two subschemes in $\Pn$
which corresponds to $[\N_0]$ via $\N$-type resolution with $\N_0 \neq 0$.
By theorem 1.20, $\L$ has a minimal element $X_0$. If we apply corollary
3.10 with $X=X_0$, we see that the admissible characters $\s \geq \g_{X_0}$
index all the cohomology preserving deformation classes of $\L$.
By propositions 2.6 and 2.8, they can also be indexed by functions
$\eta$ and $\theta$ which satisfy certain conditions. \par
\definition{Definition 3.12}
Let $\L$ be an even linkage class of subschemes of codimension two in $\Pn$.
Assume that $\L$ corresponds to $[\N_0]$ via $\N$-type resolutions, with
$\N_0 \neq 0$. Let $X_0$ be a minimal element for $\L$. Then we define
$s_0(\L)=s_0(X_0), s_1(\L)=s_1(X_0)$, $t_1(\L)=t_1(X_0)$ and $e(\L)=e(X_0)$.
\enddefinition
\definition{Remark 3.13}
We note that the integers above are well-defined.
Indeed, if $X_0, X_0^\prime$ are minimal elements for $\L$, then theorem 1.20
shows that $X_0$ and $X_0^\prime$ have the same cohomology, and hence
$s_0(X_0)=s_0(X_0^\prime)$, $s_1(X_0)=s_1(X_0^\prime)$ and $e(X_0)=
e(X_0^\prime)$. To see that $t_1(\L)$ is well-defined, we apply
theorem 1.21 to see that any minimal subscheme $X_0$ links to a minimal
subscheme $Y_0$ for the dual linkage class
via a hypersurface $S$ of degree $s_0(\L)$ and another hypersurface $T$.
Since the degrees of $X_0$ and $Y_0$ are determined by the even linkage class
$\L$, so is the degree of $T$. This degree is $t_1(\L)$.
\enddefinition
\definition{Definition 3.14}
Let $\L, X_0$ be as in definition 3.12. Let $X \in \L$. We define the
{\it height of $X$} to be the unique nonnegative integer $h_X$ such that
$M_X^i \cong M_{X_0}^i(-h_X)$ for each $1 \leq i \leq n-2$.
\enddefinition
\definition{Definition 3.15}
Let $\L, X_0$ be as in definition 3.12.
Let $\g_0=\g_{X_0}$. Then we define $\eta_X=\eta_{\g_0,\g_X,h_X}$, as defined
in proposition 2.6. Similarly we define $\theta_X=\theta_{\g_0,\g_X,h_X}$, as
defined in proposition 2.9.
\enddefinition
\proclaim{Proposition 3.16}
Let $\L$ be an even linkage class as in definition 3.12. If $X \in \L$,
then \roster
\item $s_0(X)=s_0(\L)+\sum{\theta_X(l)}$
\item $s_1(X)=\text{min}\{s_1(\L)+h_X, (\theta_X)_a\}$
\item $e(X)=\text{max}\{e(\L)+h_X, (\eta_X)_o-n\}$
\endroster
\endproclaim
\demo{Proof}
Suppose that $X$ is obtained from $X^\prime$ by a basic double link of
type $(s,1)$. Corollary 3.7 shows that $\sum{\theta_{X^\prime,X}(l)}=1$
if $s > s_0(X^\prime)$, and is equal to $0$ if $s=s_0(X^\prime)$. Noting
that $e(X)=\text{max}\{l:H^{n-1}(\I_X(l)) \neq 0\}$ and looking at the
way an $\N$-type resolution changes in going from $X^\prime$ to $X$
(proposition 1.14), we find that $e(X)=\text{max}\{e(X^\prime)+1,s-n\}$.
Using these two results, we find that if we obtain $X$ from a minimal
element $X_0$ by a sequence of basic double links of height one, then
the first and third formulas above hold. By proposition 3.9 and corollary
3.10, we can find ${\overline X}$ which has the same cohomology as $X$ and
is obtained from a minimal element by a sequence of basic double links of
height one, which proves (1) and (3). \par
For the second relationship, we combine the formulas for $\eta$
(proposition 2.6) and $\theta$ (proposition 2.9) to get
$$\g_X(l)=\g_{X_0}(l)-{l \choose 0}+{{l-h_X} \choose 0}+
{{l-s_0(X)} \choose 0}-{{l-s_0(X_0)-h_X} \choose 0}+\theta(l).$$
If we subtract the $\theta(l)$ part, we can use the definitions of
$s_0$ and $s_1$ to see that what remains is equal to $0$ for
$s_0(X) \leq l < s_1(X_0)+h_X$ and positive for $l=s_1(X_0)+h_X$.
The second equality follows.
\enddemo
\definition{Remark 3.17}
Another invariant is described in \cite{2} which also describes the
location of a subscheme $X$ in $\L$ with respect to a minimal element $X_0$.
This invariant is a sequence of integers $\{b,g_2,g_3,\dots,g_r\}$ where
$b \geq 0, g_i \leq g_{i+1}$ and $r+b=h_X-1$. For the reader who is familiar
with this invariant, we briefly describe how to go back and forth via
the invariant $\theta_X$. \par
If we know $\{b,g_2,g_3,\dots,g_r\}$, then we obtain the function $\theta_X$
by the rule
$$\theta_X(l)=\#\{k:g_k+r-k=l\}$$
In this case we have $\sum{\theta(l)}=r-1$ and $s_0(X)=s_0(X_0)+r-1$. \par
Conversely, given the function $\theta_X$ with $m=\sum{\theta_X(l)}$, we
compute the corresponding invariant $\{b,g_2,g_3,\dots,g_r\}$ by
$b=h_X-m-1, r=m-1$ and
$$g_k=r-k+max\{s:\sum_{l \geq s}{\theta(l)}>r-k\}.$$
\enddefinition
\proclaim{Lemma 3.18}
Let $\L$ be as in definition 3.12 and let $X \in \L$ be a subscheme
which links to $Y$ by a complete intersection of hypersurfaces of degrees
$s$ and $t$. Letting $s_0=s_0(\L)$ and $t_1=t_1(\L)$, we have that
$$h_X+h_Y=s+t-s_0-t_1.$$
\endproclaim
\demo{Proof}
Let $X_0$ and $Y_0$ be a pair of minimal elements for their even linkage
classes, which are linked by hypersurfaces of degrees $s_0$ and $t_1$ (see
remark 3.13). If $X_0$ has an $\N$-type resolution giving an exact sequence
$$0 @>>> \P_0 @>>> \N_0 @>>> \I_{X_0} @>>> 0$$
then in applying proposition 1.6 we get an $\E$-type resolution for $Y$ which
gives an exact sequence
$$0 @>>> \N_0^\vee(-s_0-t_1) @>>> \P_0^\vee(-s_0-t_1) \oplus \O(-s_0) \oplus
\O(-t_1) @>>> \I_{Y_0} @>>> 0.$$
$X$ has an $\N$-type resolution which involves $\N_0(-h_X)$ and $Y$ has an
$\E$-type resolution which involves $\N_0^\vee(-s_0-t_1-h_Y)$. Since
$X$ and $Y$ are linked by hypersurfaces of degrees $s$ and $t$, we
can apply proposition 1.6 to get another $\E$-type resolution for $Y$
which involves $\N_0^\vee(h_X-s-t)$. By looking at the higher Rao
modules of $Y$ (which are not all zero because we assumed that $\L$ was
not the ACM class), we see that the twists must agree, so we have
$h_X-s-t=-s_0-t_1-h_Y$.
\enddemo
\proclaim{Theorem 3.19}
Let $\L, X_0$ be as in definition 3.12. Let $s_0=s_0(\L)$ and $t_1=t_1(\L)$ and
let $X \in \L$ be a subscheme which links to $Y$ by a complete intersection
of hypersurfaces of degrees $s$ and $t$.
Then the functions $\eta_X$ and $\eta_Y$ are related by the formula
$$\eta_X(l)-\eta_Y(s+t-1-l)=$$
$${{l-s} \choose 0}+{{l-t} \choose 0}-
{{l-s_0-h_X} \choose 0}-{{l-t_1-h_X} \choose 0}.$$
\endproclaim
\demo{Proof}
By theorem 1.21 there exists a minimal element $Y_0$ for the dual linkage
class to $\L$ which links to $X_0$ by surfaces of degrees $s_0$ and $t_1$.
Conceptually, the situation is represented by a diagram
$$\CD
X @>{s,t}>> Y \\
@AAA  @AAA \\
X_0 @>{s_0,t_1}>> Y_0
\endCD$$
in which the horizontal arrows denote linkages and the vertical arrows
denote dominations. By theorem 1.20, $X_0$ has an $\N$-type resolution
which gives an exact sequence
$$0 @>>> \P_0 @>>> \N_0 @>>> \I_{X_0} @>>> 0.$$
We next apply proposition 1.6 to the $\N$-type resolution for $X_0$
to get an $\E$-type resolution for $Y_0$ which gives an exact sequence
$$0 @>>> \N_0^\vee(-s_0-t_1) @>>> \P_0^\vee(-s_0-t_1) \oplus \O(-s_0) \oplus
\O(-t_1) @>>> \I_{Y_0} @>>> 0.$$ \par
Now obtain $X$ from  $X_0$ by an even number of links. If we apply
proposition 1.6 for each link, we get a resolution for $\I_X$
$$0 @>>> \P_0(-k) \oplus \P_X @>>> \N_0(-k) \oplus \Q_X $$
where $\P_X$ and $\Q_X$ are dissoci\'e. In looking at the higher Rao modules,
we see that $k=h_Y$. Remark 2.9 shows that
$\eta_X(l) = \Delta^n\I_X(l)-\Delta^n\I_{X_0}(l-h_X)$, which can be rewritten
as as $\Delta^n\Q_X(l)-\Delta^n\P_X(l)$ in view of the two resolutions.
Similarly, we apply proposition 1.6 while doing an even number of linkages
to get from $Y_0$ to $Y$, and obtain a resolution for $\I_Y$ of the form
$$0 @>>> \N_0^\vee(-s_0-t_1-h_Y) \oplus \P_Y @>>> \P_0^\vee(-s_0-t_1-h_Y)
\oplus
\O(-s_0-h_Y) \oplus \O(-t_1-h_Y) \oplus \Q_Y$$
where $\P_Y$ and $\Q_Y$ are dissoci\'e. As above, we have
$\eta_Y(l)=\Delta^n\Q_Y(l)-\Delta^n\P_Y(l)$. \par
Now we apply proposition 1.6 to the linkage between $X$ and $Y$. After using
the formula of lemma 3.18 to simplify some expressions, we obtain
another resolution for $\I_X$ of the form
$$\CD
0 \\
@VVV \\
\P_0(-h_X) \oplus \Q_Y^\vee(-s-t) \oplus \O(-t_1-h_X)
\oplus \O(-s_0-h_X) \\
@VVV \\
\N_0(-h_X) \oplus \P_Y^\vee(-s-t) \oplus \O(-s) \oplus \O(-t) \\
\endCD$$
Comparing with the other resolution for $\I_X$, we see that
$\eta_X(l)$ is given by
$$\Delta^n\P_Y^\vee(l-s-t)-\Delta^n\Q_Y^\vee(l-s-t)+$$
$${{l-s} \choose 0}+{{l-t} \choose 0}-{{l-s_0-h_X} \choose 0}-
{{l-t_1-h_X} \choose 0}.$$
Since $\P_Y$ and $\Q_Y$ are dissoci\'e sheaves of the same rank, we can write
$$-\Delta^n\Q_Y^\vee(l-s-t)+\Delta^n\P_Y^\vee(l-s-t)=
\Delta^n\Q_Y(s+t-1-l)-\Delta^n\P_Y(s+t-1-l),$$
which is just $\eta_Y(s+t-1-l)$. Substituting this into the formula for
$\eta_X$ gives the result.
\enddemo
\proclaim{Corollary 3.20}
With the same hypotheses as theorem 3.19, assume further that $s=s_0(X)$. Then
$\theta_X$ is related to $\eta_Y$ by the formula
$$\theta_X(l)=\eta_Y(s+t-1-l)+{{l-t} \choose 0}-{{l-t_1-h_X} \choose 0}$$
\endproclaim
\demo{Proof}
This follows immediately from the formula relating $\theta_X$ and $\eta_X$ of
proposition 2.9.
\enddemo
%
%
\heading $4.$ Subschemes which satisfy $s_1=t_1$ \endheading
In this section we study the subschemes $X$ in an even linkage class $\L$
which satisfy $s_1(X)=t_1(X)$. As in the previous section, we restrict
to the case in which $\L$ is not the class of arithmetically Cohen-Macaulay
subschemes. A purely numerical criteria is found for
the existence of such a subscheme in terms of the invariant $\theta$.
This is achieved by finding a sharp lower bound for $t_1(X)$, from
which the theorem follows immediately. \par
   We mentioned earlier that we always have the inequality
$s_1(X) \leq t_1(X)$. Before going on to the results in this section, we
give an example to show that this inequality can be strict.
\definition{Example 4.1}
Let $C_0$ be the disjoint union of two lines in $\Pthree$ and let $C$ be
an elementary double link from $C_0$ of type $(8,1)$. Then one easily computes
that $s_0(C)=s_1(C)=3$, but we must have $t_1(C)>3$ because $C$ has degree
$10$ and hence cannot be contained in a complete intersection of cubic
surfaces. In fact, it can be shown that $t_1(C)=8$; $t_1(C) \leq 8$ because
$C$ was an elementary double link
of type $(8,1)$ from $C_0$. To see that $t_1(C) \geq 8$,
suppose that $C$ links to a curve $D$ via surfaces of degrees $3$ and $t$.
Then the degree of $D$ is $3t-10$, while the height of $D$ is $t-2$. Since
$D$ is in the even linkage class of two skew lines, we see that when the
height of $D$ is $t-2$, the degree of $D$ is at least $2+2(t-2)=2t-2$.
Comparing these two estimates on the degree of $D$, we see that $t \geq 8$.
\enddefinition
\proclaim{Lemma 4.2}
Let $F_a, F_b \subset \Pn$ be hypersurfaces of degrees $a \leq b$ which
meet properly. If $F_c$ is any other hypersurface, then there exists
a hypersurface $G_b$ of degree $b$ which contains $F_a \cap F_b$ and
meets $F_c$ properly.
\endproclaim
\demo{Proof}
Let $X$ denote the complete intersection $F_a \cap F_b$ and consider
the projective space $\p=\p{H^0(\I_X(b))}$. Let
$f$ be the equation of the hypersurface $F_c$
and factor $f$ into its irreducible factors as $\Pi_{i=1}^rp_i$. Let
$[p_i] \subset \p$ denote the linear subspace consisting of all
nonzero multiples of $p_i$ (note that $[p_i]$ may be empty). Then each $[p_i]$
is a proper linear subspace of $\p$, since for
each $i$ either $p_i$ does not divide the equation for $F_b$ or $p_i$ does
not divide the equation for $F_a$, and hence does not divide some multiple
of that equation of degree $b$.
Since these are proper, the union $\cup_{i=1}^r[p_i]$ is a closed proper
subspace of $\p$. It follows that for the general element
$g \in H^0(\I_X(b))$, none of the $p_i$ divide $g$, hence the
intersection $Z(f) \cap Z(g)$ is proper.
\enddemo
\proclaim{Corollary 4.3}
Let $X$ be a codimension two subscheme of $\Pn$ which is contained in
surfaces $F_a,F_b$ of degrees $a \leq b$ which meet properly. Then
$t_1(X) \leq b$.
\endproclaim
\demo{Proof}
Apply lemma 4.2 with $F_c=S$, $S$ a surface of degree $s_0(X)$ which
contains $X$.
\enddemo
\definition{Remark 4.4}
One consequence of corollary 4.3 is that for any subscheme $X$,
the least degree of a complete intersection of two hypersurfaces
containing $X$ is precisely $s_0(X)t_1(X)$, and that any pair of
surfaces giving such a complete intersection have degrees $s_0(X)$ and
$t_1(X)$.
\enddefinition
\definition{Notation 4.5}
For the remainder of this section, we let $\L$ denote a fixed even linkage
class of codimension two subschemes in $\Pn$ such that $\L$ corresponds to
$[\N_0]$ via $\N$-type resolution with $\N_0 \neq 0$.
We set $s_0=s_0(\L), s_1=s_1(\L)$ and $t_1=t_1(\L)$.
\enddefinition
\proclaim{Theorem 4.6}
Let $X \in \L$ and let $\theta=\theta_X$. Then we have the following
sharp lower bounds on $t_1(X)$. \roster
\item If $\theta_o < t_1+h_X-1$, then $t_1(X) \geq t_1+h_X$.
\item If $\theta_o \geq t_1+h_X-1$, then
$t_1(X) \geq max\{l:\theta(l) \neq 0 \text{ and } \theta(l-1)=0\}$
\endroster
\endproclaim
\demo{Proof}
Let $s=s_0(X)$ and $t=t_1(X)$. Then there exist hypersurfaces $S,T$
of degrees $s,t$ such that $S \cap T$ links $X$ to another scheme $Y$.
Corollary 3.20 gives the formula
$$\theta(l)=\eta_Y(s+t-1-l)+{{l-t} \choose 0}-{{l-t_1-h_X} \choose 0}.$$
Now suppose that $\theta_o < t_1+h_X-1$ (This includes the case when
$\theta=0$). If $t<t_1+h_X$, then formula above shows that
$\eta_Y(s+t-t_1-h_X)=-1$, which contradicts the fact that
$\eta_Y$ is a nonnegative function. Thus $t \geq t_1+h_X$.
\par
Now consider the case when $\theta_o \geq t_1+h_X-1$. Since $\eta_Y$ is
connected in degrees $< s_0+h_Y$, $\eta_Y(s+t-1-l)$ is connected
in degrees $\geq s+t-t_1-h_Y=s_0+h_X$. Looking at the formula
relating $\theta$ and
$\eta_Y$, we conclude that $\theta$ is connected in degrees $\geq t$
(consider the cases $t \leq t_1+h_X$ and $t > t_1+h_X$ separately).
It follows that
$t \geq \text{max}\{l:\theta(l) \neq 0 \text{ and } \theta(l-1)=0\}$.
\par
To see that these bounds are sharp, we construct examples by induction on
$\sum{\theta(l)}$. If this sum is zero, then we can obtain a subscheme
$X^\prime$ by directly linking to a minimal element $Y_0$ for the dual
linkage class hypersurfaces of degrees $s_0$ and $t_1+h_X$. From corollary
4.3 it is clear that $t_1(X^\prime) \leq t_1+h_X$. \par
Now assume that $\theta \neq 0$. Let
$r=\text{max}\{l:\theta(l) \neq 0 \text{ and } \theta(l-1)=0\}$ and let
$w=\theta_o-r+1$ (this is the width of the rightmost connected piece of
$\theta$). Define $\theta^\prime$ by
$$\theta^\prime(l)=\theta(l+w)-{{l-r-w} \choose 0}+
{{l-\theta_o-1} \choose 0}.$$
Because $\theta$ satisfies the criteria of proposition 2.9 for the
height $h_X$, it is easily seen that $\theta^\prime$ does as well for
the height $h_X-w$. Hence $\theta^\prime=\theta_{X^\prime}$ for some
subscheme in $\L$ of height $h_X-w$. If we define $r^\prime$ for $X^\prime$
analogously to the definition of $r$ for $X$, we can use our induction
hypothesis to assume that
$t_1(X^\prime) \leq \text{max}\{r^\prime, t_1+h_X-w\}$. \par
Using proposition 3.6, we see that if ${\overline X}$ is obtained from
$X^\prime$
by an elementary double link of type $(r,w)$, then ${\overline X}$ is in the
cohomology preserving deformation class of $X$. If $\theta_o < t_1+h_X-1$,
then $r < t_1+h_X-w-1$ and we can first link $X^\prime$ to $Y$ by hypersurfaces
of degrees $r, t_1+h_X-w$ and then link $Y$ to ${\overline X}$ by hypersurfaces
of degrees $r, t_1+h_X$. From corollary 4.3, we see that
$t_1({\overline X}) \leq t_1+h_X$. On the other hand, if
$\theta_o \geq t_1+h_X-1$, then $r \geq r^\prime$, and we can use
hypersurfaces of degrees $r, s_0(X^\prime)$ to link $X^\prime$ to $Y$, and
then use hypersurfaces of degrees $r, s_0(X^\prime)+w$ to link $Y$ to
${\overline X}$. In this case
we have that $r \geq s_0(X^\prime)$, so corollary 4.3 shows that
$t_1({\overline X}) \leq r$. In either case, we obtain sharp bounds.
\enddemo
\proclaim{Theorem 4.7}
Let $X \in \L$. Then $X$ can be deformed with constant cohomology through
subschemes in $\L$ to ${\overline X}$ satisfying $s_1({\overline X})=
t_1({\overline X})$ if and only if $\theta_X$ is connected about the
interval $[s_1+h_X, t_1+h_X-1]$.
\endproclaim
\demo{Proof}
Let $\theta=\theta_X$. Supposing that $X \in \L$ satisfies $t_1(X)=s_1(X)$,
we first show that $\theta_X$ is connected about $[s_1+h_X, t_1+h_X-1]$.
If $\theta = 0$, then proposition 3.16 shows that $s_0(X)=s_0$ and
$s_1(X)=s_1+h_X$. If $s_1 < t_1$, then lemma 3.18 shows that $X$ links
to a subscheme of negative height, a contradiction. Hence $s_1=t_1$, the
interval in empty, and $\theta$ is connected about $[s_1+h_X,t_1+h_X-1]$,
as we wanted. \par
Now assume that $\theta \neq 0$. In this case, proposition 3.16 shows
that $s_1(X) \leq \theta_a$. If $\theta_o < t_1+h_X-1$, then from
theorem 4.6 we get
$$s_1(X) \leq \theta_a \leq \theta_o < t_1+h_X \leq t_1(X)$$
and hence $s_1(X)<t_1(X)$, contradiction. We may henceforth assume
that $\theta_o \geq t_1+h_X-1$. If $\theta$ is not connected, then
by the second case of theorem 4.6 we have
$$s_1(X) \leq \theta_a < \text{max}\{l:\theta(l) \neq 0 \text{ and }
\theta(l-1)=0\} \leq t_1(X),$$
which also contradicts our assumption, so $\theta$ is connected.
Finally, if $\theta_a > s_1+h_X$, then
$s_1(X)=s_1+h_X < \theta_a \leq t_1(X)$, which also contradicts assumption.
We have shown that $\theta_o \geq t_1+h_X-1$, $\theta_a \leq s_1+h_X$, and
$\theta$ is connected. This proves the forward part of the theorem.
Now suppose that $\theta$ is connected about $[s_1+h_X, t_1+h_X-1]$. We
will use the sharpness of the bounds in theorem 4.6 to show that $X$
deforms to ${\overline X}$ satisfying $s_1({\overline X})=
t_1({\overline X})$. If $\theta_o < t_1-h_X-1$, then $\theta$ can only
be connected about $[s_1+h_X, t_1+h_X-1]$ if $\theta=0$ and the
interval is empty. In this case, $s_1=t_1$, and using the sharpness
of theorem 4.6, we can find ${\overline X}$ in the deformation class
of $X$ with $t_1({\overline X})=t_1+h_X=s_1+h_X=s_1(X)$ (this last
equality follows from proposition 3.16 and the fact that $\theta=0$).
\par
Now assume that $\theta_o > t_1+h_X-1$ (in particular, $\theta \neq 0$).
The condition on $\theta$ shows that $\theta_a \leq s_1+h_X$, hence
$s_1(X)=\theta_a$ by proposition 3.16. Using the sharpness of the second
case of theorem 4.6, we can find ${\overline X}$ in the deformation class
of $X$ such that $t_1({\overline X})=\theta_a=s_1({\overline X})$. This
finishes the proof.
\enddemo
\proclaim{Corollary 4.8}
With the hypotheses of theorem 4.7, assume that $X$ links to $Y$ by
hypersurfaces of degrees $s=s_0(X)$ and $t=t_1(X)$. If $\theta=\theta_X$ is
the zero function, then so is $\theta_Y$. If $\theta \neq 0$, then
$$\theta(l)=\theta_Y(s+t-1-l)+{{l-\theta_a} \choose 0}-
{{l-\theta_o-1} \choose 0}$$
\endproclaim
\demo{Proof}
If $\theta = 0$, then the formula of corollary 3.20 makes it clear that
$\eta_Y=0$ hence $\theta_Y=0$, so we go on to the case where $\theta \neq 0$.
In this case, theorem 4.7 shows that $t=s_1(X)=\theta_a$. Corollary 3.20
gives us the formula
$$\theta(l)=\eta_Y(s+t-1-l)+{{l-\theta_a} \choose 0}-{{l-t_1-h_X} \choose 0}.$$
Now we consider two cases. If $(\eta_Y)_a \geq s_0+h_Y$, then
$\eta_Y=\theta_Y$, $\eta_Y(s+t-1-l)=0$ for $l \geq t_1+h_X$ (by lemma 3.18),
and hence $\theta_o \leq t_1+h_X-1$. The property of theorem 4.6 shows that
$\theta_o = t_1+h_X-1$, and we deduce the result. \par
Now we consider the case $(\eta_Y)_a < s_0+h_Y$. In this case $\eta_Y(s+t-1-l)$
is nonzero for some $l \geq t_1+h_X$ and we must have
$(\eta_Y)_a=s+t-1-\theta_o$. Since $\eta_Y=\theta_Y+1$ on the interval
$[(\eta_Y)_a,s_0+h_Y-1]$ and $t_1+h_X=s+t-(s_0+h_Y)$, we can write
$$\theta_Y(s+t-1-l)=\theta_Y(s+t-1-l)+{{l-t_1-h_X} \choose 0}-{{l-s-t-\eta_a}
\choose 0}.$$
Substituting in the formula for $(\eta_Y)_a$ (in terms of $\theta_o$) into
the above expression gives the formula.
\enddemo
\proclaim{Corollary 4.9}
Let $\M \subset \L$ denote the subset of subschemes $X$ such that
$s_1(X)=t_1(X)$. Then $\M$ has a minimal element with respect to domination.
\endproclaim
\demo{Proof}
We produce a minimal element for $\M$ as follows. Let $X_0$ be an (absolute)
minimal element for $\L$. If we link $X_0$ to a minimal element $Y_0$ of
the dual class via hypersurfaces $S,T$ of degrees $s_0, t_1$ (possible
by theorem 1.21) and then link $Y_0$ to a subscheme $X_1$ via hypersurfaces
$S^\prime, T$ of degrees $s_0+t_1-s_1, t_1$, the subscheme $X_1$ is
obtained from $X_0$ by an elementary double link of type $(t_1, t_1-s_1)$.
{}From corollary 4.3 it is evident that $t_1(X_1) \leq t_1$. On the other
hand, we also find from proposition 3.6 that
$$\eta_{X_1}(l)={{l-t_1} \choose 0}-{{l-2t_1+s_1} \choose 0}.$$
Since $(\eta_{X_1})_a=t_1 \geq s_0+h_{X_1}$, we have $\theta_X=\eta_X$
and proposition 3.16 shows that $s_1(X_1)=t_1$. Thus $s_1(X_1)=t_1(X_1)$
and $X_1 \in \M$. \par
Let $X \in \M$. Since $s_1(X)=t_1(X)$, theorem 4.7 shows that $\theta_X$
is connected about $[s_1+h_X,t_1+h_X-1]$. Since $\theta_{X_1}$ is the
identity function for the interval $[s_1+h_{X_1}, t_1+h_{X_1}-1]$, we
easily check that the function
$$\theta(l)=\theta_X(l)-\theta_{X_1}(l-h_X+h_{X_1})$$
is nonnegative. Since $\sum{\theta_{X_1}(l)}=h_{X_1}$, it is immediate
that $\sum{\theta(l)} \leq h_X-h_{X_1}$. We can now apply proposition 2.10
to see that $X_1 \leq_{h_X-h_{X_1}} X$. Hence each element in $\M$
dominates $X_1$, so $X_1$ is a minimal element for $\M$. \par
\enddemo
\definition{Remark 4.10}
The reader might wonder if the subposet $\M \subset \L$ of proposition 4.9
satisfies the Lazarsfeld-Rao property. The answer is no.
For example, if $X_1$ is the minimal
element of $\M$, let $X_2$ be obtained from $X_1$ by a basic
double link of type $(t_1+h_{X_1}+1,1)$. In this case $\theta_{X_2}$ is
connected about $[s_1+h_{X_2}, t_1+h_{X_2}-1]$, but it is impossible
for $X_2$ to satisfy $s_1(X_2)=t_1(X_2)$ because all hypersurfaces of
degrees $< t_1+h_{X_1}+1$ which contain $X_2$ share a common hyperplane.
\enddefinition
%
%
\heading $5.$ Integral Subschemes in an Even Linkage Class \endheading
In this section we give a necessary conditions for subschemes to be
integral and investigate to what extent these conditions are sufficient.
For this section we work in a fixed even linkage class $\L$ of codimension
two subschemes of $\Pn$ which corresponds to $[\N_0]$ with $\N_0 \neq 0$.
To simplify notation we set $s_0=s_0(\L), s_1=s_1(\L), t_1=t_1(\L)$ and
$e=e(\L)$. Before working in $\L$, we prove a few technical lemmas to
allow for more conceptual proofs of the main results.
\par
\proclaim{Proposition 5.1}
Let $X \subset \Pn$ be an integral subscheme of codimension two. Let
$t \in \Z$ be an integer such that $t=s_0(X)$ or $t \geq s_1(X)$. Then
the general hypersurface $H$ containing $X$ of degree $t$ enjoys the
following properties. \roster
\item $H$ is integral.
\item $X$ is not contained in the singular locus of $H$.
\item $X$ is generically Cartier on $H$.
\endroster
\endproclaim
\demo{Proof}
Let $H \supset X$ be a hypersurface of degree $t = s_0(X)$ whose equation
is $h$. If $h=fg$ is reducible, then $X \subset Z(f)$ or
$X \subset Z(g)$ because $X$ is integral, which would contradict the
minimality of the degree of $H$. We conclude that $H$ is integral.
Now suppose that $X$ is contained in the singular locus of $H$.
$X$ is reduced, so it would lie on the hypersurfaces with equations
$\partial{h}/\partial{X_i}$, where $X_i$ are the homogeneous coordinates.
These equations must all be zero, as otherwise the minimality of the
degree of $H$ is contradicted. From Euler's formula
$\sum(\partial{h}/\partial{X_i})X_i = th$, we see that $t=0$,
the characteristic of $k$ is $p>0$, and $h$ is a $p^{th}$ power of another
equation $g$. Since $h$ is irreducible, we again reach a contradiction, and
conclude that $X$ is not contained in the singular locus of $H$. Thus we have
proved the first two conditions in the special case when $t=s_0(X)$. \par
Now fix a hypersurface $S \supset X$ of degree $s_0(X)$ and let $T$ be a
hypersurface of degree $s_1(X)$ which is not a multiple of $S$. If the equation
of $T$ is $f$, we again find that $f$ is not reducible, for if $f=gh$,
then $X \subset Z(g)$ or $X \subset Z(h)$ would give a hypersurface of
degree $< s_1(X)$ which meets $S$ properly. Thus $T$ is integral and $X$ is
contained in the complete intersection $S \cap T$.
Because $\I_{S \cap T}(s_1(X))$ is generated
by global sections, the general hypersurface $H \supset S \cap T$
of degree $> s_1(X)$ is integral. Since integrality
among hypersurfaces of a fixed degree is an open condition, the general
hypersurface $H \supset X$ of degree $> s_1(X)$ is integral. This proves
the first property. \par
The set of $h \in \p{H^0(\I_X(t))}$ such that
$X \subset Z(\partial{h}/\partial{X_i})$ for each $i$ is closed. On the
other hand, this closed set must be proper, since if $S \supset X$ is
a hypersurface of minimal degree, then the general union of $S$ with
a nonsingular hypersurface which meets $X$ properly will give a hypersurface
whose singular locus does not contain $X$. This checks the second property.
If $x \in X \cap T_{reg}$, then the ideal of $X$ at $x$ is a height one
prime ideal in the regular local ring $\O_{T,x}$, hence is principal by
Krull's Hauptidealsatz. This shows that $X$ is generically Cartier on
each hypersurface satisfying the second property.
\enddemo
\proclaim{Proposition 5.2}
Let $X \in \L$ be an integral subscheme with $0 \neq \theta_X=\theta$.
Then $\theta_a \leq s_0+h_X$.
\endproclaim
\demo{Proof}
By way of contradiction, assume that $X$ is integral and $\theta_a > s_0+h_X$.
Applying proposition 5.1, we see that $X$ is linked to a subscheme $Y$
by integral hypersurfaces $S,T$ of degrees $s=s_0(X),t=s_1(X)$.
Corollary 4.8 shows that $\sum{\theta_Y(l)} < \sum{\theta(l)}$, hence
$s_0(Y) < s_0(X)=s$. Since $S$ is an integral hypersurface of degree $s$,
we must have that $t_1(Y) \leq s$. By theorem 4.7, we have that
$t=\theta_a > s_0+h_X$, hence $s < t_1+h_Y$ by lemma 3.18, and
$t_1(Y) < t_1+h_Y$. On the other hand, corollary 4.8 shows that
$\theta_Y(s+t-1-l)=0$ for $l < \theta_a$. In particular, this holds for
$l \leq s_0+h_X$, so $\theta_Y(l)=0$ for
$l \geq s+t-1-s_0-h_X=t_1+h_Y-1$ and $(\theta_Y)_o < t_1+h_Y-1$.
theorem 4.6 now gives that $t_1(Y) \geq t_1+h_Y$, a contradiction.
\enddemo
Now we prove another necessary condition on integral subschemes, but
this one requires more work. In \cite{9}, Lazarsfeld and Rao show
that if $e(C)+4 < s_0(C)$ for a curve $C$ in $\Pthree$, then $C$ is
the unique minimal curve in its even linkage class. I will show that
if $X \in \L$ and $e+n+1+h_X < s_0(X)$, then $\L$ has a unique minimal
element $X_0$, and that $X_0 \subset X$. First I must recall a few
results on moduli for subschemes in fixed even linkage classes. \par
\proclaim{Proposition 5.3}
Let $X \in \L$ be a subscheme. Then the subset
$$H_X=\{X^\prime:X \text{ deforms to} X^\prime \text{ through schemes in } \L
\text{ with constant cohomology }\}$$
of the Hilbert scheme is irreducible.
\endproclaim
\demo{Proof}
Fixing $X$, we can find $N$ such that $\I_X$ is $N$-regular. By the theorem
of Castelnuovo-Mumford \cite{17,p.99}, the total ideal of $X$ is generated
by its homogeneous parts of degree $\leq N$. Letting
$\Q = \oplus_{l \leq N}{\O(-l)^{h^0(\I_X(l))}}$, we can find an $\E$-type
resolution for $X$ of the form
$$ 0 @>>> \E @>>> \Q @>>> \I_X @>>> 0$$
and each $X^\prime \in H_X$ has an $\E$-type resolution of the same form.
Let $V$ be the parameter space of all morphisms $\{\E @>\varphi>> \Q\}$.
$V$ is a smooth projective variety and comes equipped with a universal
morphism $p^*\E @>\varphi>> p^*\Q$, where $p:V \times \Pn @>>> \Pn$ is
the second projection. \par
Let $W \subset \Pn$ denote the subset where $\E$ is locally free. By
\cite{6, corollary 1.4} the complement of $W$ is of codimension $\geq 3$.
On the open set $V \times W \subset V \times \Pn$ we have the complex
$$0 @>>> p^*\E @>\varphi>> p^*\Q @>{\Lambda^r\varphi^\vee \otimes 1}>>
p^*\O @>>> 0.$$
By \cite{17, corollary 4.6}, the set $U$ consisting of $v \in V$ such that
$0 @>>> p^*\E \otimes k(v) @>>> p^*\Q \otimes k(v) @>>> p^*\O$ is
exact is open. Let $j:U \times W \subset U \times \Pn$ denote the
inclusion map. Applying $j_*$ extends the complex to $U \times \Pn$.
Let $\X$ be the subscheme defined by the cokernel of the
last map in the complex. For each $u \in U$, we have a sequence
$$0 @>>> \E @>{\varphi_u}>> \Q @>>> \O @>>> \O_{\X_u} @>>> 0,$$
which is exact on $W$. Since $\E, \Q,$ and $\O$ are reflexive, and the
codimension of $W$ is $\geq 3$, this sequence is exact on $\Pn$.
In particular, the Hilbert polynomial of the fibres is constant, and
$\X$ is a flat family over $U$. Hence there is an
induced map $U @>>> \H_X$, where $\H_X$ is the Hilbert scheme for subschemes
in $\Pn$ with the same Hilbert polynomial as $X$. We take $H_X$ to be the
image with induced reduced structure.
\enddemo
\proclaim{Proposition 5.4}
If $U \subset H_X$ is an open set such that $X_0 \subset X^\prime$ for
each $X^\prime \in U$, then $X_0 \subset X^\prime$ for all
$X^\prime \in H_X$.
\endproclaim
\demo{Proof}
Since $H_X$ is irreducible by proposition 5.3, it suffices to show that
the set of $\{X^\prime \in H_X: X_0 \subset X^\prime\}$ is closed.
For this, we will show that the corresponding set is closed in $\H_X$, the
Hilbert scheme for all subschemes of $\Pn$ having the same Hilbert polynomial
as $X$. If $\F$ is the flag scheme for all inclusions $X_0 \subset X$, then
we have natural projections
$$\CD
\F @>p>> \H_{X_0} \\
@VqVV @. @. \\
\H_X @. @. \\
\endCD $$
The inverse image of $\{X_0\}$ under $p$ is closed in $\F$, and the image
of this closed set under $q$ is closed in $\H_X$, because $q$ is projective.
It follows that the set of interest is closed in $\H_X$, and hence so is
its intersection with $H_X$.
\enddemo
\proclaim{Lemma 5.5}
Let $X \subset \Pn \times H @>\pi>> H$ be a family of subschemes with
constant cohomology from an even linkage class $\L$ of codimension two
subschemes in $\Pn$. Then the function $t_1(X_h)$ is upper
semicontinuous on $H$.
\endproclaim
\demo{Proof}
Since it suffices to show this on each irreducible component of $H$,
we can base extend by the irreducible components with reduced induced
structure to reduce to the case where $H$ is integral. To show that
$t_1(X_h)$ is upper semicontinuous, we must show that for each $h \in H$,
there is an open neighborhood $h \in U \subset H$ such that
$t_1(X_k) \leq t_1(X_h)$ for each $k \in U$. \par
Fix $h \in H$, and let $t=t_1(X_h)$. Since the dimension $h^0(\I_{X_k}(t))$
is constant for $k \in H$, we see by Grauert's theorem
\cite{6, III, corollary 12.9} that the sheaf
$\F=\pi_*(\I_X \otimes q^*(\O(t)))$ is locally free on $H$. Letting
${\p}\F @>f>> H$ be the corresponding vector bundle over $H$, and base
extending the universal family by $f$, we obtain a flag scheme for
inclusions $X_k \subset T$, where $k \in H$ and $T$ is a hypersurface
of degree $t$. We can similarly find a bundle ${\p}\G @>g>> H$ which
parametrizes inclusions $X_k \subset S$ with $k \in H$ and $S$ a hypersurface
of degree $s$. \par
   By definition of $t_1(X_h)$, we can find hypersurfaces $S,T$ containing
$X$ of degrees $s,t$ such that $S \cap T$ has dimension $n-2$. Since
${\p}\F @>f>> H$ is a vector bundle, we can find an open set $U_t \subset H$
and a local section $\s_t:U_t @>>> f^{-1}U_t$ such that $\s_t(h)$ corresponds
to the inclusion $X_h \subset T$. In this way we obtain a flat family
$$\CD
{\tilde X_t} \subset {\tilde T} @>>> \Pn \times U_t \\
@. @VVV \\
@. U_t \\
\endCD $$
We can obtain a similar family
${\tilde X_s} \subset {\tilde S} \subset \Pn \times U_s @>>> U_s$ whose
fibre at $h$ corresponds to the inclusion $X_h \subset S$. \par
Letting $U = U_s \cap U_t$ be the intersection of the open sets, we
obtain families ${\tilde X} \subset {\tilde T}$ and
${\tilde X} \subset {\tilde S}$ over $U$. The map
$$ {\tilde S} \cap {\tilde T} @>>> U$$
is surjective and the fibre over $h$ has dimension $n-2$. Now we can use
the semicontinuity of the dimension of the fibres of this morphism
\cite{6, II, exercise 3.22} to find an open set $h \in V \subset U$
over which the fibres have dimension $\leq n-2$. Clearly the intersections
${\tilde S_k} \cap {\tilde T_k}$ have dimension $\geq n-2$, so for
$k \in V$ this intersection has dimension exactly $n-2$, which shows
that $t_1(X_k) \leq t$.
\enddemo
\proclaim{Proposition 5.6}
Let $\L$ be an even linkage class such that $0 < \delta=s_0-e-n-1$.
Then $\L$ has a unique minimal element $X_0$. Further, If $Y \in \L^\vee$
satisfies $(\eta_Y)_o < t_1+h_Y+\delta$ and $Y$ links to $X$ by hypersurfaces
of degrees $< t_1+h_Y+\delta$, then $X_0 \subset X$.
\endproclaim
\demo{Proof}
   For the first statement, we use the same proof as in \cite{9}. One
can use a minimal $\E$-type resolution for a dual curve $Y_0$. Upon
linking to $X_0$, applying proposition 1.6 to get a resolution for
$\I_{X_0}$ and cancelling the two summands corresponding to the
hypersurfaces used for the linkage, we obtain an $\N$-type resolution
$$0 @>>> \P @>>> \N @>>> \I_{X_0} @>>> 0$$
in which we can write $\P \cong \oplus{\O(l)^{p(l)}}$, with $p(l)=0$
for $l > e+n+1$. In this situation, any other injection $\P @>>> \N$
must have the same image, hence the isomorphism class of the ideal
sheaf $\I_{X_0}$ is uniquely determined. Since $X_0 \subset \Pn$ is
of codimension two, it follows that $X_0$ is unique. \par
For the statement about $Y$, we proceed by induction on $h_Y$. If
$h_Y=0$, then $Y=Y_0$ is a minimal curve for $\L^\vee$. We know
that $Y_0$ links directly to $X_0$ by hypersurfaces $S,T$ of degrees
$s_0, t_1$ by theorem 1.21. It follows from \cite{7} that there is an
exact sequence
$$0 @>>> \I_{S \cap T} @>>> \I_{Y_0} @>>> \omega_{X_0}(n+1-s_0-t_1) @>>> 0.$$
Twisting by any $l < t_1+\delta$ and taking global sections gives an
isomorphism $H^0(\I_{S \cap T}(l)) \cong H^0(\I_{Y_0}(l))$, which
shows that $H^0(\I_{Y_0}(l))$ is generated by $S$ and $T$. It follows
that if $F_a \cap F_b \supset Y_0$ is a complete intersection with
$a,b < t_1+\delta$, then $F_a \cap F_b \supset S \cap T = X_0 \cup Y_0$,
which proves the proposition for the case $h_Y=0$. \par
Now suppose that $h_Y > 0$. We can find $Y_1$ of height $h_Y-1$ for
which an elementary double link of type $(A,1)$ will be in the cohomology
preserving deformation class of $Y$ by proposition 3.9.
Since $(\eta_Y)_o < t_1+h_Y+\delta$, it is clear that
$(\eta_{Y_1})_o < t_1+h_{Y_1}+\delta$, so our induction hypothesis
applies to $Y_1$. Using the sharpness of the bounds in theorem 4.6,
we can find such a $Y_1$ with $t_1(Y_1) < t_1+h_{Y_1}+\delta$, Hence
we can link $Y_1$ to $X_1$ by hypersurfaces of degrees
$A,B=t_1+h_{Y_1}+\delta-1$. By induction hypothesis we have
that $X_0 \subset X_1$. By lemma 5.5, there is an open subset in $H_{X_1}$
consisting of subschemes $X^\prime$ which link by hypersurfaces of degrees
$A$ and $B$ to subschemes $Y^\prime \in H_{Y_1}$, hence we
see that the general $X \in H_{X_1}$ contains $X_0$. By proposition 5.5, we
have that all $X^\prime \in H_{X_1}$ must contain $X_0$. \par
Now we return to $Y$. Suppose that $Y$ links to some $X$ by
hypersurfaces of degrees $a,b$, which
are both less than $t_1+h_Y+\delta$. Then $t_1(Y)<t_1+h_Y+\delta$ by
corollary 4.3. This shows that $Y$ links to subschemes $X^\prime$ by
hypersurfaces of degrees $A,B+1$. Since $Y$ was in the deformation class
of an elementary double link from $Y_1$ of type $(A,1)$, we see that
$X^\prime \in H_{X_1}$, hence $X_0 \subset X^\prime$. It follows that
the general hypersurface $T$ of degree $B+1$ containing $Y$ also contains
$Y \cup X_0$, and hence this is true for all such $T$ as this condition
is closed. In particular, $X$ must contain $X_0$, for otherwise we could
find a hypersurface $F_a$ of degree $a$ which does not contain $Y \cup X_0$,
and the union of $F_a$ with a general hypersurface of degree
$t_1+h_Y+\delta-1-a$ would also not contain $Y \cup X_0$, a contradiction.
\enddemo
\proclaim{Proposition 5.7}
Let $X \in \L$. If $s_0(X) > e+n+1+h_X$, then \roster
\item $\L$ has a unique minimal element $X_0$
\item $X_0 \subset X$
\endroster
\endproclaim
\demo{Proof}
Since $s_0(X)=s_0+\sum{\theta_X(l)} \leq s_0+h_X$, we have that
$s_0 > e+n+1$, so $\L$ has a unique minimal element $X_0$ by
proposition 5.7. Set $\delta=s_0-e-n-1$. Letting $s=s_0(X)$ and
choosing $t> \text{max}\{t_1(X), t_1+h_X\}$, we can link $X$ to $Y$
by hypersurfaces of degrees $s$ and $t$. The hypothesis on $t$ and
the formula of corollary 3.20 show that $\eta_Y(l)=0$ for $l \geq t$,
so $(\eta_Y)_o < t$. Further, the hypothesis $s > e+n+1+h_X$ combined
with lemma 3.18 show that $t < t_1+h_Y+\delta$, so we can apply
proposition 5.7 to $Y$ with $a=s$ and $b=t$ to conclude that $X_0 \subset X$.
\enddemo
\proclaim{Theorem 5.8}
Let $X \in \L$ be integral and assume that $X$ is not minimal. Then \roster
\item $\theta_X$ is connected about $[s_0+h_X, t_1+h_X-1]$
\item $s_0(X) \leq e+n+1+h_X$.
\endroster
\endproclaim
\demo{Proof}
The first condition is a consequence of theorem 4.7 and proposition 5.2.
If $s_0(X) > e+n+1+h_X$, then proposition 5.7 shows that $\L$ contains
a minimal element $X_0$, which is contained in $X$. This contradicts the
assumption that $X$ is integral, unless $X$ is minimal.
\enddemo
\proclaim{Proposition 5.9}
Let $T \subset \Pn$ be an integral hypersurface, $Y \subset T$ a
generically Cartier generalized divisor. Let $m$ be an integer such
that the linear system ${\p}H^0(\I_{Y,T}(m))$ cuts out a scheme
$Y^\prime \supset Y$ which differs from $Y$ on a set of codimension $>1$ in
$T$ and that $H^0(\I_{Y,T}(m-1)) \neq 0$.
Then the general complete intersection of $T$ with a hypersurface
$Z(f) \supset Y$ of degree $m$ links $Y$ geometrically to an integral
subscheme $X$.
\endproclaim
\demo{Proof}
The linear system $V=H^0(\I_{Y,T}(m))$ has $Y^\prime$ as base locus, hence
defines a morphism $T-Y^\prime @>\s>> {\p}V$. By hypothesis there exists
$0 \neq s \in H^0(\I_{Y,T}(m-1))$. Let $S \subset \Pn$ be a hypersurface
whose restriction to $T$ is the scheme of zeros of $s$ and let
$W \subset V$ denote the sublinear system of multiples
of $s$. $W$ gives an embedding of $T-S$ into ${\p}W$, which can be factored
as $T-S \subset \Pn-S @>>> {\p}W$. The first map is a
closed immersion and the second is an open immersion, hence the composite
map is unramified. This composite map can also be factored
$T-S @>>> T-Y^\prime @>\s>> {\p}V @>\pi>> {\p}W$ where $\pi$
is a projection from a linear subspace. It follows that $\s$ is unramified
when restricted to $T-S$ and that the dimension of the image of $\s$ is at
least two. \par
We are now in position to apply Jouanolou's Bertini theorem
\cite{8, theorem 6.10}. Let $H_f$ be a general hyperplane in ${\p}V$ which
corresponds to $f \in H^0(\I_{Y,T}(m))$. By Jouanolou's Bertini theorem,
$\s^{-1}(H_f)$ is geometrically irreducible and $\s^{-1}(H_f)-S$ is
reduced. Also, since $Y^\prime$ is generically Cartier on $T$ and
$\I_{Y^\prime,T}(m)$ is generated by its global sections, the general
$f \in H^0(\I_{Y,T}(m))$ generates $I_{Y^\prime,T}$ at its generic points of
codimension one in $T$ and meets $S-Y^\prime$ in codimension $>1$. Putting
these
facts together, we find that for general $f$, $Z(f) \cap T=Y \cup X$
where $X$ has no common component with $S$, is reduced and geometrically
irreducible when restricted to $T-S$. In other words, $X$ links to $Y$
geometrically and is integral.
\enddemo
\proclaim{Proposition 5.10}
Let $X \in \L$ be an integral and $h > 0$ an integer. Let $t \leq e(X)+n+1+h$
be an integer such that $t=s_0(C)$ or $t \geq t_1(C)$. Then the general
elementary double link $X^\prime$ of type $(t,h)$ is integral.
\endproclaim
\demo{Proof}
Let $t$ be as in the theorem, and let $T$ be a general surface
of degree $t$ which contains $X$. Then $T$ is integral and $X$ is generically
Cartier on $T$ by proposition 5.1.
By proposition 5.10, a general hypersurface $H$ of large degree $d$
links $X$ geometrically to an integral subscheme $Y$.
We have an isomorphism $\I_Y/{\I_{S \cap T}} \cong \omega_X(4-d-t)$.
Restricting the ideal sheaves to the surface $T$ yields the exact sequence
$$0 \ra \I_{S \cap T,T} \ra \I_{Y,T} \ra \omega_X(4-d-t) \ra 0$$
Note that $\I_{S \cap T,T} \cong \O_T(-d)$ and that the map on the left
is just multiplication by the equation for $S$. If we twist this sequence
by $d+h$ we get an exact sequence on global sections
$$0 \ra H^0(\O_T(h)) \ra H^0(\I_{Y,T}(d+h)) \ra H^0(\omega_X(n+1-t+h)) \ra 0.$$
Since $t \leq e(X)+n+1+h$, the last cohomology group is nonzero, hence not
every element of $V=H^0(\I_{Y,T}(d+h))$ is a multiple of the equation for
$S$, although all multiples of this equation by linear forms is contained
in $V$. It follows that the map $V \otimes \O_T(-d-h) \ra \O_T$ defines
the ideal sheaf of a curve $Y^\prime$ such that
$Y \subset Y^\prime \subset Y \cup X$. Further, the second inclusion is
proper because not every element of $V$ was a multiple of the equation of
$S$. Since $X$ is integral, the proper subscheme $Y^\prime \cap X$ has
codimension $>2$. \par
Since we also have $H^0(\I_{D,T}(d)) \neq 0$, the hypotheses of proposition
5.10 hold and we find that $Y$ links geometrically via $T$ and a hypersurface
$S^\prime$ of degree $d+h$ to a subscheme $X^\prime$ which is integral.
This completes the proof.
\enddemo
\proclaim{Theorem 5.11}
Let $X,Y \in \L$ such that $X$ is integral, $X \leq Y$, and $Y$ satisfies
the conclusion of theorem 5.8. Then $Y$ can be deformed with constant
cohomlogy through schemes in $\L$ to ${\overline Y}$ with ${\overline Y}$
integral.
\endproclaim
\demo{Proof}
We induct on the relative height $h=h_Y-h_X$.
Let $\eta=\eta_{X,Y}, \theta=\theta_{X,Y}$.
Define $r=\text{min}\{l:\eta(l) \neq 0 \text{ and } \eta(l+1)=0\}$ and let
$A=r-h_Y+h_X+1, w=r-\eta_a+1$. If $X^\prime$ is obtained from $X$ by an
elementary double link of type $(A,w)$, then applying propositions 2.7 and
3.6 shows that $X^\prime \leq Y$. Since $h_Y-h_{X^\prime}=h-w<h$, we can
apply the induction hypothesis once we show that $X^\prime$ can be taken
integral. For this we will check the conditions of proposition 5.11. \par
First we check that $A=s_0(X)$ or $A \geq s_1(X)$. If $A < s_0(X)$, then
because $s_0(X) \leq (\eta_X)_a$ we find that $r < (\eta_X)_a+h_Y-h_X-1$.
Proposition 2.7 gives us the formula $\eta_Y(l)=\eta_X(l+h_Y-h_X)+\eta(l)$,
so the inequality would imply that $\eta_Y$ is not connected, a contradiction,
so we must have $A \geq s_0(X)$. \par
Now assume that $A > s_0(X)$. Then $r > s_0(X)+h_Y-h_X-1$, which implies
that $\theta(r) \neq 0$. Since $\eta(r+1)=0$, we also have $\theta(r+1)=0$.
Now consider the formula of proposition 2.10.
$$\theta_Y(l)=\theta_X(l+h_Y-h_X)+\theta(l)$$
If $\theta_X \neq 0$, then we must have $(\theta_X)_a+h_Y-h_X \leq r+1$ for
$\theta_Y$ to be connected. In this case, theorem 4.7 shows that
$s_1(X) \leq r-h_Y+h_X+1$, which is equivalent to $A \geq s_1(X)$.
If $\theta_X =0$, then since $\theta_Y$ must be connected about
$[s_0+h_Y, t_1+h_Y-1]$, we find that $r \geq t_1+h_Y-1 \geq s_1+h_Y-1$,
which shows that $A \geq s_1+h_X$, which is $\geq s_1(X)$ by proposition
2.11. We conclude that $A=s_0(X)$ or $A \geq s_1(X)$. \par
Secondly we check that $A \leq e(X)+n+1+w$. In view of the formula
$e(X)=\text{max}\{e+h_X, (\eta_X)_o-n\}$, we find that this is equivalent
to showing that $\eta_a \leq \text{max}\{e+n+1, (\eta_X)_o+1-h_X\}+h_Y$.
If $\eta_X=0$, then the fact that $\eta_Y$ is connected about
$[e+n+1+h_Y, t_1+h_Y-1]$ gives this fact. If $\eta_X \neq 0$, the condition
on $\eta_X$ from theorem 5.9 gives that $(\eta_X)_o \geq t_1+h_X-1$.
Since $t_1 \geq e+n+1$, we need to check that
$\eta_a \leq (\eta_X)_o+h_Y-h_X+1$.
If this is not the case, then the formula which relates $\eta, \eta_X$ and
$\eta_Y$ shows that $\eta_Y$ is not connected, contradicting the conditions
of theorem 5.9. Hence $A \leq e(X)+n+1+w$ and our proof is complete.
\enddemo
\definition{Example 5.12}
For any even linkage class which has an integral minimal element, the
conditions of theorem 5.8 are both necessary and sufficient for the
existence of an integral subscheme. For example, Lazarsfeld and Rao show
that if $C$ is a smooth curve, then the general embedding of $C$ in
$\Pthree$ of large degree will have image which is minimal in its even
linkage class.
\enddefinition
\definition{Example 5.13}
Let $\L$ be the even linkage class of $4$ skew lines which lie on a
quadric surface $Q$ in $\Pthree$. One can calculate that $s_0=2,
s_1=t_1=4$, and $e=-1$. In this case, it turns out that there are
two minimal integral curves. One can be obtained from the $4$ lines
by taking a general elementary double link of type $(2,1)$. These
curves are of type $(1,5)$ on the smooth quadric, and hence have smooth
connected representatives. The other can be obtained by a general
elementary double link of type $(4,2)$. This can be produced by first
linking the $4$ skew lines to another set of $4$ skew lines via $Q$ and
a quartic, and then the second set of $4$ skew lines can be geometrically
linked to an integral curve by two quartic surfaces via proposition 5.9.
In this case, the conditions of theorem 5.8 are both necessary and
sufficient.
\enddefinition
\definition{Example 5.14}
Using the construction of Hartshorne and Hirschowitz, one can show that
the general union of $10$ lines meeting at $9$ points generizes to a
smooth rational curve of degree $10$. In this case, both the rational
curve and the union of the $10$ lines have seminatural cohomology, and
we can compute that $e=-1$ and $s_0=5$, so that both of these curves are
unique minimal elements in their even linkage classes by proposition 5.6.
On the other hand, one even linkage class has a minimal integral curve and
the other doesn't. We conclude that it is not possible to get purely
cohomological criteria for the existence of integral subschemes in an
even linkage class. In the case of the rational curve, there is just
one class with a minimal integral curve of height $h>0$. It can be
obtained from the rational curve by an elementary link of type $(5,2)$.
\enddefinition
\definition{Example 5.15}
Consider the even linkage class of a double line $C_0$ with arithmetic genus
$g = -4$. In \cite{17, theorem 8.2.7}, it is shown that the even linkage
class of $C_0$ has two minimal integral curves, which can be obtained from
$C_0$ as follows. The first can be obtained by first linking $C_0$ to $D_1$
by quadric surfaces, and then linking $D_0$ to $C_1$ by a cubic and quintic
surface. In general, the curve $C_2$ produced in this way will be smooth and
connected. The second minimal integral curve can be obtained by first linking
$C_0$ to $D_2$ by a quadric and quartic surface, and then linking $D_2$ to
$C_2$ by a cubic and a sextic surface. $C_2$ can again be taken to be smooth
and connected. It is also shown that not every deformation class containing
an integral curve contains a smooth connected curve, so we cannot hope
to get results like theorem 5.8 and theorem 5.11 for smooth connected
curves without adding some extra hypotheses.
\enddefinition

\enddocument